\documentclass[11pt]{article}
\usepackage{amscd,amsmath, amssymb, fancyhdr}

\newcommand{\version}{version 4.0,\ \  Jan 22, 2011}

\numberwithin{equation}{section}

\def\eqref#1{(\ref{#1})}
\newcommand{\goth}{\mathfrak}

\newcommand{\arrow}{{\:\longrightarrow\:}}
\newcommand{\Z}{{\Bbb Z}}
\newcommand{\C}{{\Bbb C}}
\newcommand{\R}{{\Bbb R}}
\newcommand{\Q}{{\Bbb Q}}
\renewcommand{\H}{{\Bbb H}}
\newcommand{\6}{\partial}
\def\1{\sqrt{-1}\:}
\newcommand{\restrict}[1]{{\left|_{{\phantom{|}\!\!}_{#1}}\right.}}
\newcommand{\cntrct}                
{\hspace{2pt}\raisebox{1pt}{\text{$\lrcorner$}}\hspace{2pt}}

\newcommand{\calo}{{\cal O}}

\renewcommand{\bar}{\overline}
\renewcommand{\phi}{\varphi}
\renewcommand{\epsilon}{\varepsilon}

\newcommand{\im}{\operatorname{im}}

\newcounter{Mycounter}[section]
\newcounter{lemma}[section]
\renewcommand{\thelemma}{{Lemma \thesection.\arabic{lemma}}}
\newcommand{\lemma}{%
     \setcounter{lemma}{\value{Mycounter}}
     \refstepcounter{lemma}
     \stepcounter{Mycounter}
     {\noindent \bf \thelemma:\ }}

\newcounter{claim}[section]
\renewcommand{\theclaim}{{Claim \thesection.\arabic{claim}}}
\newcommand{\claim}{%
     \setcounter{claim}{\value{Mycounter}}
     \refstepcounter{claim}
     \stepcounter{Mycounter}
     {\noindent \bf \theclaim:\ }}

\newcounter{sublemma}[section]
\newcounter{corollary}[section]
\renewcommand{\thecorollary}{{Corollary \thesection.\arabic{corollary}}}
\newcommand{\corollary}{%
     \setcounter{corollary}{\value{Mycounter}}
     \refstepcounter{corollary}
     \stepcounter{Mycounter}
     {\noindent \bf \thecorollary:\ }}

\newcounter{theorem}[section]
\renewcommand{\thetheorem}{{Theorem \thesection.\arabic{theorem}}}
\newcommand{\theorem}{%
     \setcounter{theorem}{\value{Mycounter}}
     \refstepcounter{theorem}
     \stepcounter{Mycounter}
     {\noindent \bf \thetheorem:\ }}

\newcounter{conjecture}[section]
\newcounter{proposition}[section]
\renewcommand{\theproposition}
       {{Proposition \thesection.\arabic{proposition}}}
\newcommand{\proposition}{%
     \setcounter{proposition}{\value{Mycounter}}
     \refstepcounter{proposition}
     \stepcounter{Mycounter}
     {\noindent \bf \theproposition:\ }}

\newcounter{definition}[section]
\renewcommand{\thedefinition}
       {{Definition~\thesection.\arabic{definition}}}
\newcommand{\definition}{%
     \setcounter{definition}{\value{Mycounter}}
     \refstepcounter{definition}
     \stepcounter{Mycounter}
     {\noindent \bf \thedefinition:\ }}

\newcounter{example}[section]
\newcounter{remark}[section]
\renewcommand{\theremark}{{Remark \thesection.\arabic{remark}}}
\newcommand{\remark}{%
     \setcounter{remark}{\value{Mycounter}}
     \refstepcounter{remark}
     \stepcounter{Mycounter}
     {\noindent \bf \theremark:\ }}

\newcounter{problem}[section]
\newcounter{question}[section]
\renewcommand{\thequestion}{{Question \thesection.\arabic{question}}}
\newcommand{\question}{%
     \setcounter{question}{\value{Mycounter}}
     \refstepcounter{question}
     \stepcounter{Mycounter}
     {\noindent \bf \thequestion:\ }}

\makeatletter

\@addtoreset{equation}{section}
\@addtoreset{footnote}{section}
\makeatother

\def\blacksquare{\hbox{\vrule width 5pt height 5pt depth 0pt}}
\def\endproof{\blacksquare}

\begin{document}
\begin{center}
{\LARGE\bf Oeljeklaus-Toma manifolds admitting no complex subvarieties
}\\

\hfill

Liviu Ornea\footnote{Partially supported by a PN2-IDEI grant, nr. 525, and by
Tokyo Institute of Technology during summer of 2010.}
and Misha Verbitsky\footnote{Partially supported by the
 RFBR grant 10-01-93113-NCNIL-a,
RFBR grant 09-01-00242-a, Science Foundation of 
the SU-HSE award No. 10-09-0015 and AG Laboratory HSE, 
RF government grant, ag. 11.G34.31.0023\\

{\bf Keywords:} Locally
conformally K\"ahler manifold,
K\"ahler potential, positive bundle, complex subvariety, Inoue surface.

{\bf 2000 Mathematics Subject
Classification:} { 53C55.}}

\hfill

{\em To Professor Vasile Br\^{\i}nz\u anescu at his sixty-fifth birthday}

\end{center}

\hfill

{\small
\hspace{0.15\linewidth}
\begin{minipage}[t]{0.7\linewidth}
{\bf Abstract} \\
The Oeljeklaus-Toma (OT-) manifolds are complex manifolds
constructed by Oeljeklaus and Toma from certain 
number fields, and generalizing the Inoue surfaces $S_m$.
On each OT-manifold we construct a holomorphic line 
bundle with semipositive curvature form $\omega_0$
and trivial  Chern class. Using this form, we prove that the OT-manifolds
admitting a locally conformally K\"ahler structure
have no non-trivial complex subvarieties.
The proof is based on the Strong Approximation theorem
for number fields, which implies that any leaf of the
null-foliation of $\omega_0$ is Zariski dense.
\end{minipage}
}

\tableofcontents


\section{Introduction}


\subsection{OT-manifolds and their subvarieties}

The Oeljeklaus-Toma (OT-) manifolds are an important class
of compact complex manifolds not admitting a K\"ahler 
metric. They were discovered by 
Oeljeklaus and Toma in  2005 (\cite{_Oeljeklaus_Toma_}).
The construction of OT-manifolds uses 
the Dirichlet unit theorem from number theory
(Subsection \ref{_OT_constru_Subsection_};
see \cite{_Parton-Vuletescu_} for additional details
of this construction and many related questions).
Starting from a degree 3 number field, one obtains
a 2-dimensional OT-manifold known as Inoue surface $S_m$ 
(see \cite{_Inoue_}).

For some number fields, the OT-manifolds
are locally conformally K\"ahler. A locally conformally K\"ahler
(LCK) structure on a complex manifold is a K\"ahler metric
on its universal cover $\tilde M$, such that the deck transform
maps act on $\tilde M$ by homotheties.
The OT-manifolds serve an important function
in the theory of LCK manifolds, providing a counterexample
to a longstanding conjecture of I. Vaisman, \cite{vai}, who asked whether
there exists a compact, non-K\"ahler LCK-manifold $M$ 
with all odd Betti numbers even: $b_{2p+1}(M)\equiv 0 \mod 2$.
The Oeljeklaus-Toma manifolds in dimension $3$ are the only known examples of
compact LCK-manifolds with even odd Betti numbers, $b_1=b_5=2, b_3=0$. 

An OT-manifold is LCK if  it is
constructed from a number field $K$
which has precisely 2 complex (non-real) embeddings,
that is, two distinct 
homomorphisms $K \stackrel {\sigma, \bar\sigma} \arrow \C$. If the OT manifold
has at least 4 complex embeddings and exactly one real, then it is not LCK. The
remaining case is not yet decided.

Inoue surfaces $S_m$ have no curves. We give a generalization
of this theorem, proving that an OT-manifold which is
locally conformally K\"ahler has no non-trivial
complex subvarieties. In particular, it has no non-constant meromorphic
functions (as a meromorphic function without polar set would be holomorphic, and
hence constant).

\hfill

\question Is there any OT (non-LCK) manifold that has non-constant meromorphic
functions?

\hfill

The idea of the proof of this result is quite simple.
We construct a holomorphic Hermitian line bundle, called
{\bf the weight bundle}, on any OT-manifold $M$. This bundle
is topologically trivial, and has semipositive curvature
form $\omega_0$. The weight bundle also admits a flat 
connection, compatible with the holomorphic structure. 

To learn about complex subvarieties of an OT-manifold,
we study the zero-foliation $\Sigma$ of $\omega_0$, proving that 
all its leaves are Zariski dense in $M$. For an
OT-manifold $M$ constructed from a number field $K$
admitting exactly $2t$ distinct complex (non-real)
embeddings to $\C$, the leaves of $\Sigma$ are
$t$-dimensional. When $t=1$, $M$ is locally conformally
K\"ahler, and $\Sigma$ is one-dimensional. In this case, we 
prove that for any positive-dimensional complex subvariety 
$Z\subset M$, $Z$ contains with each point $z\in Z$ a leaf 
$\Sigma_z$ passing through $z$. Since all leaves of
$\Sigma$ are Zariski dense, the same is true for $Z$.

The weight bundle $L$ is quite useful for many
other purposes. As it was done in
\cite{_Verbitsky:Sta_Elli_}, one can take the $\alpha$-th 
tensor power of $L$, denoted by $L^\alpha$, for any real $\alpha$;
this power is well defined, because $L$ is equipped with a
natural $C^\infty$-trivialization. The Gauduchon degree
$\deg_g$ of $L^\alpha$, taken with respect to any
Gauduchon metric, satisfies 
$\frac 1 \alpha \deg_g L^\alpha= \deg_g L>0$, hence
$M$ admits a line bundle with any prescribed Gauduchon
degree. This implies, in particular, that the connected
component of the Picard group $\mathop{Pic}(M)$ is non-compact. 
Also, this implies that any vector bundle on $M$
has degree zero after tensoring with an appropriate power
of $L$; this is useful for the study of Hermitian-Einstein
bundles on $M$, providing useful tools
for the classification of stable bundles, and, eventually,
coherent sheaves on $M$.

A similar argument was used in \cite{_Verbitsky:Sta_Elli_}
to study holomorphic vector bundles and subvarieties
on homogeneous elliptic fibrations, such as Calabi-Eckmann 
manifolds and quasi-regular Vaisman manifolds. 
We pose two questions,
very much unsolved, but quite natural in the context
presented by \cite{_Verbitsky:Sta_Elli_} and the present paper.
Notice that from their construction it is clear that
OT-manifolds are affine flat, that is, equipped with a
flat, affine, torsion-free connection.

It is shown in \cite[Remark 1.7]{_Oeljeklaus_Toma_} that some OT manifolds
admit 
a holomorphic foliation with compact leaves which are
again OT manifolds. Hence, it is natural to 
pose the following:

\hfill

\question
Are the ones described in \cite[Remark 1.7]{_Oeljeklaus_Toma_} the only OT
manifolds with compact 
complex subvarieties? Can we classify these subvarieties?
Are they always totally geodesic with respect to the
flat affine connection?

\hfill

\question
Does there exist a stable holomorphic vector bundle
of rank $>1$ on any OT-manifold of dimension $>2$? Do all holomorphic
vector bundles admit a flat connection, compatible
with the holomorphic structure?

\hfill

\remark
It is well known that generic complex tori have no 
non-trivial complex subvarieties. In \cite{_Verbitsky_torus_}, 
it was shown that all stable bundles
on a generic complex torus of dimension $>2$ have rank 1, and
all holomorphic vector bundles admit flat connections. As for compact complex
surface of non-K\"ahlerian type, it is proven in \cite{vuli} that  
stable holomorphic $2$-bundles with $c_1=0$ and $c_2=n$ exist for any $n>0$.

\subsection{Number theory and the construction of OT-manifolds}
\label{_OT_constru_Subsection_}


Let $[K:\Q]$ be a number field, that is, a finite extension
of $\Q$, of degree $n$, with $\sigma_1, ..., \sigma_s$
the real embeddings of $K$ into $\C$, and $\sigma_{s+1}, ..., \sigma_n$
the complex embeddings. Since the
complex embeddings of $K$ into $\C$
occur in pairs of complex conjugate embeddings, the number $n-s$ is
even, $n-s=2t$.  Let
$\sigma=(\sigma_1,\ldots,\sigma_n):K\rightarrow
\C^{s+t}$ be the corresponding group homomorphism. 

Let $\calo_K$ be the ring of algebraic integers of $K$,
$\calo_K^*$ its multiplicative group of units and
$\calo^{*,+}_K$ the group of units which are positive in all
the real embeddings of $K$. 

Denote by $\H$ the upper complex half-plane. 
Using the Dirichlet's unit theorem, 
Oeljeklaus and Toma proved that
$\calo_K\rtimes\calo^{*,+}_K$ acts freely on ${\Bbb H}^s\times
\C^t$ by 
\begin{equation*}
\begin{split}
T_a(z_i)&=(z_i+\sigma_i(a)), \quad i=1,\ldots,s+2t,\quad a\in \calo_K,\\
R_u(z_i)&=(\sigma_i(u)z_i), \quad i=1,\ldots,s+2t,\quad u\in \calo_K^{*,+}.
\end{split}
\end{equation*}
(see \cite{_Oeljeklaus_Toma_}, \cite{_Parton-Vuletescu_}).
Moreover, an {\em admissible}  subgroup $U\subset
\calo^{*,+}_K$ can always be found such that the action of
$\Gamma:=\calo_K\rtimes U$ is also properly
discontinuous. For $t=1$, every $U$ of finite index in
$\calo^{*,+}_K$ has this property.

\hfill

\definition The manifold $M_K:=({\Bbb H}^s\times \C^t)/\Gamma$
is called an {\bf Oeljeklaus-Toma manifold}. It is a 
compact complex manifold of dimension $s+2t$.

\hfill

For $s=t=1$, $M_K$ reduces to an Inoue surface $S_m$
(where $m$ is a matrix in $\mathrm{SL}(3,\Z)$), see
\cite{_Inoue_}. The corresponding number field
$K$ is $\Q[T]/P(t)$, where $P_m(t)$
is the characteristic polynomial of the
matrix $m$. It is shown in \cite{_Oeljeklaus_Toma_} that
the manifolds $M_K$ are never K\"ahler, but that 
for $t=1$, $M_K$ is a locally conformally K\"ahler (LCK)
manifold (see \cite{drag} and the more recent survey
\cite{OV_surv} for definitions and results in LCK
geometry). We briefly explain the construction of this LCK
metric. 

Clearly, the function $\psi(z)=\prod_{i=1}^s(\im z_i)
+|z_{s+1}|^2$ is plurisubharmonic on ${\Bbb H}^s\times \C$.
It defines the  K\"ahler form $\Omega:=\6\bar\6\ \psi$
on ${\Bbb H}^s\times \C$. The group $\Gamma$ 
acts on $({\Bbb H}^s\times \C, \Omega)$ by homotheties: 
\begin{equation*}
\begin{split}
T_a^*\Omega&=\Omega,\\
R_u^*\Omega&=|\sigma_{s+1}(u)|^2\Omega.
\end{split}
\end{equation*}
  Let now $\chi:\Gamma\rightarrow \R^{>0}$ be the character
 $\chi(\gamma)=\frac{\gamma^*\Omega}{\Omega}$. We call
 {\bf automorphic} any $p$-form $\eta\in
 \Lambda^p({\Bbb H}^s\times \C)$ which satisfies
 $\gamma^*\eta=\chi(\gamma)\eta$. For any
 automorphic function $\phi$ on ${\Bbb H}^s\times \C$,   the
 quotient $\frac{\Omega}{\phi}$ is $\Gamma$-invariant and
 hence projects to an LCK metric $\omega$ on $M_K$. This form
 satisfies the equation $d\omega=\theta\wedge\omega$, for
 the closed $1$-form $\theta$ (called {\bf the Lee form})
 which is the projection on $M_K$ of
 $\tilde\theta=-d\log\phi$:
\[ 
  d\omega = -\frac {d\phi}{\phi^2} \wedge \tilde \omega =
  -d(\log \phi)\wedge \omega.
\]
It is easily seen that the function $\phi=
\prod_{i=1}^s(\im z_i)^{-1}$ is automorphic, and hence it
produces a LCK metric on $M_K$ as described above. This
LCK metric generalizes the one constructed by Tricerri on
$S_m$, \cite{_Tricerri_}.


\hfill


The main result of this paper shows that, just as Inoue
surfaces $S_m$ have no complex curves, OT-manifolds have
no complex subvarieties:

\hfill

\theorem
Let $[K:\Q]$ be a number field of degree $n=s+2$, with
$s$ real embeddings and 2 complex embeddings, and 
$M_K$ the corresponding LCK OT-manifold. Then $M_K$ has
no non-trivial complex subvarieties.

{\bf Proof:} See \ref{_no_subva_Theorem_}. \endproof

\hfill

\corollary The LCK OT-manifold $M_K$ has no non-constant meromorphic functions. 


\section{The weight bundle of an OT-manifold}


\definition
Let  $[K:\Q]$ be a number field of degree $n=s+2t$, with
$s$ real embeddings and $2t$ complex embeddings, and 
$M_K={\Bbb H}^s \times \C^t/\Gamma$ the associated OT-manifold.
Denote by $z_1, ..., z_s$ the standard complex coordinates
on ${\Bbb H}^s$, and let $\tilde\theta:= -d\log \prod_{i=1}^s(\im
z_i)$. It is easy to see that the form
$\tilde \theta$ is $\Gamma$-invariant. Therefore
it is obtained as a lift of a form $\theta$, called
{\bf the Lee form} of the OT-manifold. When $t=1$,
this is the Lee form constructed above.

\hfill

Let $M_K$ be an OT-manifold, and $\theta$ its Lee
form. Consider a trivial Hermitian line 
bundle $L$ with connection $\nabla:= \nabla_0 + \1\theta^c$,
where $\theta^c:= I\theta$, and $\nabla_0$ is the trivial
connection on $L$. Clearly, $\nabla$ is Hermitian, and
$\nabla^{0,1}= \bar\6 + \theta^{0,1}$, where $\theta^{0,1}$
is the (0,1)-part of $\theta$. 

\hfill

\claim\label{_omega_0_formula_Claim_}
In these assumptions, the curvature $\omega_0$ of $\nabla$ is
$-\1 d\theta^c$. Moreover, this form is of type (1,1).

\hfill

{\bf Proof:}  A simple computation shows that in the
standard coordinates $z_1, ... z_s, z_{s+1}, ... z_{s+t}$, $\omega_0$
can be written as follows:
\[
\omega_0 = \1\6\bar\6\log \phi = 
\1 \sum_{i=1}^s\frac{dz_i\wedge d\bar z_i} {|\im z_i|^2},
\]
\endproof

\hfill

\definition\label{weight}
Let $M_K$ be an OT-manifold, and $L$ the holomorphic
Hermitian bundle defined above. Then $L$ is called 
{\bf the weight bundle} of $M_K$. 

\hfill

We restate \ref{_omega_0_formula_Claim_} as

\hfill

\theorem\label{_omega_0_defi_Theorem_}
Let $M_K$ be an OT-manifold, and $L$ its weight bundle
with the holomorphic Hermitian structure  and the Chern connection
$\nabla$ defined above.
Consider the form $\omega_0:= \1 \nabla^2$. Then
$\omega_0$ is a semi-positive form, which can be written
in the standard coordinates $z_1, ... z_s, z_{s+1},..., z_{s+t}$
as follows:
\[
\omega_0 = \1\6\bar\6\log \phi = 
\1 \sum_{i=1}^s\frac{dz_i\wedge d\bar z_i} {|\im z_i|^2}
\]
\endproof

\hfill

\remark 
The Vaisman manifolds are, by definition, LCK manifolds
$(M,I,g)$ satisfying the additional condition
$\nabla^g\theta=0$, where $\nabla^g$ is the
Levi-Civita connection of an LCK metric $g$. 
For all Vaisman manifolds, the 2-form
$\omega_0=d\theta^c$ is semi-positive, being zero only on
the direction of $\theta^\sharp-I\theta^\sharp$. This is a
general fact, proven in \cite{_Verbitsky_vanishing_},
independent of the particular form of
$\theta$. OT-manifolds are far from being 
Vaisman (they never admit any Vaisman metric), but the particular
expression of their Lee form gives $\omega_0$ the same
property as for Vaisman manifold. This is what inspired
our construction. 

\hfill

\remark
An object of interest in conformal geometry and,
in particular, LCK geometry is the {\bf weight bundle}. It
is the real line bundle $L\arrow M$ associated to the
representation $\mathrm{GL}(2n,\R)\ni A\mapsto |\det
A|^{\frac{1}{n}}$ (see \cite{OV_surv}). Then $L$ can be
complexified and endowed with the Chern connection
$\nabla_0 + \1\theta^c$ (where $\nabla_0$ is the trivial
connection). It can be verified that
$\omega_0=\sqrt{-1}\nabla^2$, and hence $\omega_0$ can be
seen as the curvature form of this Chern connection. 
When $t=1$ and $M_K$ is the corresponding LCK OT-manifold, this construction
gives the weight bundle introduced in \ref{weight}.

\hfill

\remark
For any OT-manifold
$M$, in addition to the Chern connection $\nabla=\nabla_0 + \1\theta^c$,
the weight bundle $L$ also admits the  connection $\nabla_0 + \theta$,
which is  flat because $d\theta=0$. 
It is clear that the $(0,1)$-part of $\nabla$ 
 coincides with the $(0,1)$-part of
this flat connection.

\hfill

The following claim is obvious from the explicit form of
$\omega_0$ (\ref{_omega_0_defi_Theorem_}). 

\hfill

\claim\label{_Sigma_defi_Claim_}
In the assumptions of \ref{_omega_0_defi_Theorem_}, 
let $\tilde \Sigma$ be the holomorphic foliation on 
the covering $\tilde M_K = {\Bbb H}^s \times \C^t$
generated by the vector fields $\frac {\6}{\6z_{s+1}}, ...,\frac {\6}{\6z_{s+t}}
$.
Then:
\begin{description}
\item[(i)] The foliation $\tilde \Sigma$ is $\Gamma$-invariant, hence
it is obtained as the pullback of a holomorphic foliation
$\Sigma$ on $M_K=\tilde M_K/\Gamma$. 
\item[(ii)] The foliation $\Sigma$ is the null-space
of the form $\omega_0$ constructed above.
\end{description}

 \endproof

\hfill

\claim\label{_tangent_to_Sigma_Claim_}
Let  $[K:\Q]$ be a number field of degree $n=s+2$, with
$s$ real embeddings and $2$ complex embeddings,
$M_K$ the corresponding LCK OT-manifold, and
$\Sigma\subset TM_K$ the holomorphic foliation
defined in \ref{_Sigma_defi_Claim_}. Consider a complex closed 
subvariety $Z\subset M_K$. Then $\Sigma$ is 
tangent to $Z$ at any point of $Z$:
\begin{equation}\label{_tange_sigma_Equation_}
\forall z\in Z, \ \  \Sigma\restrict z \subset T_z Z.
\end{equation}

{\bf Proof:}   The form $\omega_0$ has
$(n-1)$ positive eigenvalues, where $n = \dim_\C M_K$,
and its zero eigenspace at $z$ is $\Sigma\restrict z$. 
Unless \eqref{_tange_sigma_Equation_} holds at $z\in Z$,
the restriction $\omega_0\restrict Z$ has $m=\dim Z$
positive eigenvalues at $z$. Then $\int_Z \omega_0^m>0$.
This is impossible, because $\omega_0$ is exact.
\endproof

\hfill

\corollary\label{_leaf_Corollary_}
In assumption of \ref{_tangent_to_Sigma_Claim_},
let $\Sigma_z$ be a leaf of $\Sigma$ passing through
$z\in Z$. Then $\Sigma_z\subset Z$.

\endproof


\section{Complex subvarieties in LCK OT-manifold}


Using \ref{_leaf_Corollary_}, we can easily 
prove the main result of this paper.

\hfill

\theorem\label{_no_subva_Theorem_}
Let $[K:\Q]$ be a number field of degree $n=s+2$, with
$s$ real embeddings and 2 complex embeddings, and let 
$M_K$ be the corresponding OT-manifold. Then $M_K$ has
no non-trivial complex subvarieties.

\hfill

{\bf Proof:} \ref {_no_subva_Theorem_} follows
from \ref{_leaf_Corollary_} and the following more general 
proposition.

\hfill

\proposition\label{_leaf_dense_Proposition_}
Let $[K:\Q]$ be a number field of degree $n=s+2t$, $t>0$, with
$s$ real embeddings and $2t$ complex embeddings, and let 
$M_K={\Bbb H}^s \times \C^t/\Gamma$ be the associated (non-K\"ahler)
OT-manifold. Let 
$\Sigma\subset TM_K$ be the foliation defined
in \ref{_Sigma_defi_Claim_}. Consider a  
leaf of $\Sigma$, and let $Z$ be its closure. Then
\begin{description}
\item[(i)] The preimage $\pi^{-1}(Z)$ of $Z$ to $\tilde M_K={\Bbb H}^s \times
\C^t$ contains the set
\[
Z_{\alpha_1, ..., \alpha_s}:=
 \{ (z_1, ..., z_s, z_{s+1},\ldots,z_{s+t})\ \ |\ \ \im z_i=\alpha_i \}
\]
for some positive numbers $(\alpha_1,\ldots, \alpha_s)\in \R^s$.
\item[(ii)] Any complex subvariety of $M_K$ containing $Z$
must coincide with $M_K$.
\end{description}

\hfill

{\bf Proof:} The implication (i) $\Rightarrow$ (ii)
is clear, because any complex
manifold containing $Z_{\alpha_1,\ldots, \alpha_s}$
must have the same dimension as $M_K$. The proof
of (i) is a bit more elaborate.

Let $\calo$ be the ring of integers in $K$.
By construction, the group $\Gamma=\pi_1(M_K)$
is a cross-product of the additive group $\calo^+$ of $\calo$
with a subgroup of the multiplicative group $\calo^*$.
Let $\tilde \Sigma$ be the pullback of the foliation
$\Sigma$ to $\tilde M_K={\Bbb H}^s \times \C^t$. A leaf of
$\tilde \Sigma$ is given as 
\[
T_{t_1, ..., t_s}:=
 \{ (z_1, ..., z_s, z_{s+1}, ..., z_{s+t})\ \ |\ \ z_i=t_i \}
\]
for some $(t_1, ..., t_s)\in {\Bbb H}^s$.
Let $\tilde Z:=\; \pi^{-1}(Z)$ be the preimage of the
corresponding closure of a leaf of $\Sigma$. Clearly, $\tilde Z$
is the closure of $\Gamma(T_{t_1, ..., t_s})$.
Therefore, to prove \ref{_leaf_dense_Proposition_} (i)
it is sufficient to show that the closure of 
$\Gamma(T_{t_1, ..., t_s})$ contains $Z_{\alpha_1, ..., \alpha_s}$.
In fact, even the smaller group $\calo^+\subset \Gamma$ will
suffice, as seen from the following lemma, which proves
\ref{_leaf_dense_Proposition_}.

\hfill

\lemma
Let $[K:\Q]$ be a number field of degree $n=s+2t$, $t>0$ with
$s$ real embeddings and $2t$ complex embeddings, and
$\tilde M_K:={\Bbb H}^s \times \C^t$, equipped with
the action of $\calo^+$ as in Subsection
\ref{_OT_constru_Subsection_}.
Consider the subset
\[
T_{t_1, ..., t_s}:=
 \{ (z_1, ..., z_s, z_{s+1},\ldots, z_{s+t})\ \ |\ \ z_i=t_i \}
\]
in $\tilde M_K$. Then the closure of $\calo^+(T_{t_1, ..., t_s})$
coincides with 
\[
Z_{\alpha_1, ..., \alpha_s}:=
 \{ (z_1, ..., z_s, z_{s+1},\ldots,z_{s+t})\ \ |\ \ \im z_i=\alpha_i, \}
\]
with $\alpha_i:= \im t_i$.

\hfill

{\bf Proof:}
Equivalently, we may state that the closure of 
an orbit of the standard action of $\calo^+$ in ${\Bbb H}^s$ is the set
$\{ (z_1,\ldots, z_s, z_{s+1},\ldots, z_{s+t})\ \ |\ \ \im z_i=\alpha_i \}$.
This in turn is equivalent to the following

\hfill

\lemma (cf. \cite[Claim following Lemma 2.4]{_Oeljeklaus_Toma_}) 
Let $[K:\Q]$ be a number field of degree $n=s+2t$, $t>0$ with
$s$ real embeddings $\sigma_1,\ldots, \sigma_s$
and $2t$ complex embeddings. Consider the
additive group $\calo^+$ of the corresponding ring of integers.
Let $\sigma:\; \calo^+ \arrow \R^s$ map $\xi$ to
$\sigma_1(\xi),\ldots, \sigma_s(\xi)$. Then the image
of $\calo^+$ is dense in $\R^s$.

\hfill

{\bf Proof:}\footnote%
{We are grateful to Marat Rovinsky, who kindly explained to us
this proof}
Let $K$ be a number field, $\calo_K$ its ring of integers,
${\goth P}$ the set of all prime ideals of $\calo_K$,
 $V$ the product of all
archimedean completions of $K$, and $V_1$ the product of some,
but not all, archimedean completions.
Denote by $\calo_\nu$ the completion of $\calo_K$ at $\nu\in {\goth P}$,
and let $K_\nu$ be the corresponding local field.
Consider the adele space ${\goth A}$, obtained as a subset of the  
product $V\times \prod_{\nu \in {\goth P}} K_\nu$,
where all components, except finitely many, 
belong to $\calo_\nu$, and let ${\goth A}_1$
be the image of projection of ${\goth A}$
to $V_1\times \prod_{\nu \in {\goth P}} K_\nu$.
Denote by $\tau:\; K\arrow {\goth A}_1$
the natural homomorphism, which is tautological
componentwise.

{}From the Strong Approximation theorem 
(see \cite{Kneser} or \cite[Theorem
20.4.4]{ebook}\footnote{
http://modular.fas.harvard.edu/papers/ant/html/node84.html})
it follows that the image $\tau(K)$ of $K$ 
is dense in ${\goth A}_1$. Let
\[
 \calo_{{\goth A}_1}:={\goth A}_1
\cap\left( V_1\times \prod_{\nu \in {\goth P}} \calo_\nu\right)
\]
be the set of points of ${\goth A}$,
corresponding to the integer adeles.
Clearly, $\calo_{{\goth A}_1}$ is open in ${\goth A}_1$.
Therefore, the intersection $\tau(K) \cap \calo_{{\goth A}_1}$
is dense in $\calo_{{\goth A}_1}$. On the other hand,
$\tau(K) \cap \calo_{{\goth A}_1}$ consists of those
elements of the number field which are integer
at all non-archimedean places. This gives
$\tau(K) \cap \calo_{{\goth A}_1} = \tau(\calo_K)$.
Therefore, the image of $\calo_K$
to $V_1$ is dense. 

\endproof

\hfill

\remark
The above argument actually proves that the image of
$\calo_K$ in the product $V_1$ of all archimedean completions
of $K$ except one is dense in $V_1$.

\hfill

{\bf Acknowledgements:} We are grateful to Katia Amerik for her
support. Much gratitude to Marat Rovinsky for his invaluable help in
proving the approximation lemma. Part of this work was done in
Oberwolfach during the Research in Pairs programme; we are grateful
to Oberwolfach Foundation for making it possible. 
Many thanks to Victor Vuletescu and Matei Toma for insightful 
email correspondence.

{\small

\noindent {\sc Liviu Ornea\\
University of Bucharest, Faculty of Mathematics, \\14
Academiei str., 70109 Bucharest, Romania. \emph{and}\\
Institute of Mathematics ``Simion Stoilow" of the Romanian Academy,\\
21, Calea Grivitei Street
010702-Bucharest, Romania }\\
\tt Liviu.Ornea@imar.ro, \ \ lornea@gta.math.unibuc.ro

\hfill

\noindent {\sc Misha Verbitsky\\
{\sc Laboratory of Algebraic Geometry, SU-HSE,\\
7 Vavilova Str. Moscow, Russia, 117312}\\
\tt verbit@maths.gla.ac.uk, \ \  verbit@mccme.ru\\
}}

\end{document}